
\input amstex

\documentstyle{amsppt}

\magnification=\magstep1
\hcorrection{.25in}
\advance\vsize-.75in

\define \PP{\Cal P}
\define \ZZ{\Bbb Z}

\define \st{such that}
\define \si{\sigma}
\define \sft{shift of finite type}
\define \ts{tiling system}
\define \pro{prototile}

\NoBlackBoxes

\topmatter

\title      The Symbolic Dynamics of Tiling the Integers \endtitle

\author Ethan M. Coven\\William Geller\\Sylvia Silberger\\William Thurston
\endauthor

\address Department of Mathematics,
         Wesleyan University,
         Middletown, CT  06459-0128                       \endaddress
\email   ecoven\@wesleyan.edu                              \endemail
\address Department of Mathematics,
         Indiana University-Purdue University at Indianapolis,
         Indianapolis, IN 46202-3216                          \endaddress
\email   wgeller\@math.iupui.edu 								                     \endemail
\address Department of Mathematics,
         Lafayette College,
         Easton, PA 18042                                     \endaddress
\curraddr  Department of Mathematics,
           Hofstra University,
           Hempstead, NY 11549                                \endcurraddr

\email    silbergs\@lafayette.edu                                      \endemail
\address Department of Mathematics,
         University of California Davis,
         Davis, CA 95616                                   \endaddress
\email   wpt\@math.ucdavis.edu                            \endemail

\thanks
Some of this work was done at the Mathematical Sciences Research Institute
(MSRI), where research
is supported in part by NSF grant DMS-9701755.  The first two authors thank
K.~Schmidt for useful conversations
and ideas.
\endthanks

\leftheadtext{E. COVEN, W.GELLER, S. SILBERGER, AND W. THURSTON}

\subjclass Primary 58F03 \endsubjclass

\abstract
{A finite collection~$\PP$ of finite sets {\it tiles the integers\/} iff
the integers can be expressed as a
disjoint union of translates of members of $\PP$.  We associate with such a
tiling a doubly infinite sequence with
entries from~$\PP$.  The set of all such sequences is a sofic system,
called a {\it tiling system}.  We show that,
up to powers of the shift, every shift of finite type can be realized as a
tiling system.}
\endabstract

\endtopmatter

\document

\head 1. Introduction   \endhead

For notation, terminology, and basic results of symbolic dynamics,
see the book by D.~Lind and B.~Marcus \cite{LM}.

Let $\PP = \{P_1,\dots, P_K\}$
be a finite collection of
finite subsets of the integers $\ZZ = \{\dots,-1,0,1,\dots\}$,
called {\it prototiles}.
We normalize the \pro s so that each has minimum~$0$.
A {\it tile\/} is a translate of a prototile.
A {\it tiling of the integers\/} by~$\PP$ is an expression of the integers
as a disjoint union
of tiles,
$\ZZ = \bigcup(t_j + P_{k_j})$.
Corresponding to this tiling is the point
$x = (x_i) \in \prod_{i=-\infty}^{\infty} \{1,2,\dots,K\}$
defined  by
$x_i = k$ if and only if there exists~$j$ such that
$i \in t_j + P_{k_j}$ and  $k_j = k$.
Thus we can think of a tiling as being given by a bi-infinite sequence of
colors,
where the colors are in one-to-one correspondence with the prototiles.

Let $\si$ denote the {\it shift}, $\left(\sigma(x)\right)_i = x_{i+1}$.
The collection of points corresponding to tilings of the integers by~$\PP$,
denoted~$T(\PP)$,
is closed and shift-invariant.
We call $\si: T(\PP) \to T(\PP)$ a {\it tiling system}.
We first show that every tiling system is sofic.
We then prove our main result: up to powers of the shift, every  shift of
finite type can be realized as a tiling system.

\proclaim{Main Theorem} Let $\si:\Sigma \to \Sigma$ be a  shift of finite
type.  Then there is a positive integer~$m$ and a tiling system
$\si: T \to T$ \st
 \roster
  \item $T = T_0 \cup T_1 \cup \dots \cup T_{m-1}$,
        where the ~$T_i$ are closed and cyclically permuted by the shift~$\si$.
  \item $\si^m : \Sigma \to \Sigma$ is topologically conjugate to every
$\si^m : T_i \to T_i$.
 \endroster
\endproclaim

\proclaim{Corollary}
The set of topological entropies of \ts s   is the same as that of  shifts
of finite
type, i@.e@., the set of logarithms of roots of Perron numbers.
\endproclaim

\remark{Remark}
In the sequel we will sometimes, as is common in symbolic dynamics, call
the space~$T$ a tiling system, the
space~$\Sigma$ a shift of finite type, etc.
\endremark

\head 2. Tiling systems are sofic  \endhead

Consider the following three examples.
\smallskip
\noindent (1) $\PP = \{\{0\},\{0,1\}\}$.
It is more convenient to think of~$\PP$ as~$\{R,BB\}$,  ($R$ = red, $B$ = blue).
Then $T(\PP)$ is the set of all bi-infinite indexed concatenations
of~$R$ and~$B$
such that between any two consecutive occurrences of~$R$ there is an even
number of~$B$'s,
the well-known  even system.
In this case $T(\PP)$ is also a renewal system, although we do not use that
fact here.
Recall that a {\it renewal system\/} is the collection of indexed
bi-infinite sequences
which are concatenations of a finite set of finite words from some alphabet.
In the sequel, we shall abuse notation and write $T(R,BB)$ in place of~
$T(\left\{\{0\},\{0,1\}\right\})$.

\smallskip

\noindent (2) $\PP = \{\{0\},\{0,2\}\}$, which we replace by~$\{R, B\ \_\ B\}$.
Then $T(\PP)$ is the renewal system generated by words
$R, BRB$, and~$BBBB$.

\smallskip

\noindent (3) $\PP = \{R,BB\ \_\ B,Y\_\ \_\ Y\}$,
i.e., $ \left\{\{0\},\{0,1,3\},\{0,3\}\right\}$.
In this case
$T(\PP)$ is not a renewal system.

\medskip

To  show that every tiling system is sofic, recall the proof that the even
system
$\si : T(R,BB) \to T(R,BB)$ is sofic ---
it is the image of the
\sft\ $\si : \widetilde T(R_1,B_1B_2) \to \widetilde T(R_1,B_1B_2)$
under the ``drop the subscripts'' map.
Here $\widetilde T(R_1,B_1B_2)$ is the set of all bi-infinite indexed
concatenations of ~$R_1$ and~$B_1B_2$.
We show that every ``subscripted \ts'' is a \sft.
Clearly every \ts\ can be obtained from a subscripted \ts\
by dropping the subscripts.

Formally, let $\PP = \{P_1,\dots,P_K\}$ be a finite collection of
prototiles.
 Write
$$P_k = \{0 = p_{k,1}  < p_{k,2} < \dots < p_{k,\ell_k}\}$$
and define
$\widetilde T = \widetilde T(\PP)$
on alphabet
$\{(k,\ell): 1 \le k \le K, 1 \le \ell \le \ell_k\}$,
by
$x \in \widetilde T$ iff there is a tiling of the integers by members of $\PP$,
$\ZZ = \bigcup (t_j + P_{k_j})$,
\st\ for every~$i$, there exist $j = j(i)$ and~$\ell = \ell(i)$ \st\
$i \in t_j + P_{k_j}$  and  $x_i = (k_j,\ell)$.
Equivalently, $x \in \prod \{(k,\ell): 1 \le k \le K, 1 \le \ell \le \ell_k\}$
is in~$\widetilde T$ if and only if for every~$i$, $x_i = (k,\ell)$
and $1 \le \ell' \le \ell_k$ imply
$x_{i + p_{k,\ell'} - p_{k,\ell}} = (k,\ell')$.
Informally, if $x_i$ is an element of a tile, then the other elements of that
tile appear in the appropriate places of~$x$.

The following result was proved in conversations with K.~Schmidt in Warwick
in~1994.

\proclaim{Theorem}
Every ``subscripted \ts'' is a \sft.
\endproclaim

\demo{Proof}
Let $L$ be the length  of a longest prototile in~$\PP$.  (For example, $B\
\_\ B$ has
length~$3$.)  We show that $\widetilde T = \widetilde T(\PP)$ is a \sft\
by showing that if
$x \in \prod \{(k,\ell)\}$
and every solid $L$@-word which appears in ~$x$
appears in some point of ~$\widetilde T$, then $x \in \widetilde T$.

Suppose that
every solid $L$@-word which appears in~$x$
appears in some $y \in \widetilde T$.
Let $x_i = (k,\ell)$.
Since $p_{k,\ell_k} + 1 \le L$,
there exists $y \in \widetilde T$ such that
$$    y_{i-p_{k,\ell}},\dots, y_{i-p_{k,\ell}+p_{k,\ell_k}}=
      x_{i-p_{k,\ell}},\dots, x_{i-p_{k,\ell}+p_{k,\ell_k}}.
$$
But $y \in \widetilde T$ and $y_i = (k,\ell)$, so
$y_{i-p_{k,\ell}+p_{k,\ell'}} = (k,\ell')$ for $1 \le \ell' \le \ell_k$.
Hence $x \in \widetilde T$.

Informally, suppose that every solid $L$@-word which appears in ~$x$
appears in a subscripted tiling.  Since no tile is longer than~$L$, if~$x_i$
is an element of a tile, then the other elements of that tile appear in the
appropriate places of~$x$.  Therefore $x \in  \widetilde T$.
\qed
\enddemo

\proclaim{Corollary}
Every tiling system is sofic.
\endproclaim

\remark{Remark}
We cannot use $L-1$ in the proof of the theorem.
Again let $\PP = \{R,BB\}$,
so $\widetilde T = \widetilde T(R_1,B_1B_2)$ and $L = 2$.
Every $1$@-word appearing in $x = \dots B_1B_1B_1 \dots$
appears in some point of~$\widetilde T$, but $x \notin \widetilde T$.
\endremark

\medskip

Not every sofic system can be realized as a \ts,
 as is shown by the following

\proclaim{Proposition}
A \ts\ which has a point of period~$2$ must have at least two fixed points.
\endproclaim

\demo{Proof}
The point of period~$2$  is $\dots abab \dots$,
so there are two \pro s, each of which consists entirely of even integers
or entirely of odd
integers. Both tile the even integers and hence tile the integers.
Then both $\dots aaa \dots$ and $\dots bbb \dots$ are in the \ts.
\qed
\enddemo

Similarly, if a \ts\ has a point of period~$3$ or one of period~$4$, then
it must have
at least one fixed point.  The existence of a point of period greater
than~$4$ does not
imply the existence of a fixed point.

\head 3. The main theorem
\endhead

\proclaim{Main Theorem} Let $\si:\Sigma \to \Sigma$ be a  shift of finite
type.  Then there is a positive integer~$m$ and a tiling system
$\si: T \to T$ \st
 \roster
  \item $T = T_0 \cup T_1 \cup \dots \cup T_{m-1}$,
        where the ~$T_i$ are closed and cyclically permuted by the shift~$\si$.
  \item $\si^m : \Sigma \to \Sigma$ is topologically conjugate to every
$\si^m : T_i \to T_i$.
 \endroster
\endproclaim

\demo{Proof}
We may assume that $\Sigma = \Sigma_A$, the edge shift determined by a
matrix~$A$
with nonnegative integer entries.
$A$ is the adjacency matrix of a directed graph~$G$.
The alphabet of~$\Sigma_A$ is the set of arcs (directed edges)
 of~$G$ and
$x = (x_i) \in  \Sigma_A$ if and only if
for every~$i$,
the terminal vertex of~$x_i$ is the initial vertex of ~$x_{i+1}$.

For every positive integer~$m$, $\si^m : \Sigma_A \to \Sigma_A$ is
topologically conjugate to $\si : \Sigma_{A^m} \to \Sigma_{A^m}$.
We find a tiling system
$\si : T \to T$ and a positive integer~$m$ \st\
$T = T_0 \cup T_1 \cup \dots \cup T_{m-1}$,
the $T_i$ are closed and cyclically permuted by the shift,
and $\si : \Sigma_{A^m} \to \Sigma_{A^m}$
is  topologically conjugate to $\si^m : T_0 \to T_0$.
Hence  $\si^m : \Sigma_A \to \Sigma_A$ is
topologically conjugate to every $\si^m :T_i \to T_i$.

Suppose that $A$ is $V \times V$.
Choose $n > V$ so that
$$( V\max A_{ij})^{13n} < (n+1)!$$
Let $m = 13n$. Then
every entry of~$A^m =A^{13n}$ can be written (uniquely) as
$$c_1(1!) + c_2(2!) + \dots + c_n(n!),$$
where $ 0 \le c_k \le k$ for $ 1 \le k \le n$.

 We now construct the tiling system.
The prototiles will be of two types: {\it barbells\/} and {\it racks} (to
hold barbells).
We will use the same terms for the corresponding tiles.
In the sequel we will use colors to
label prototiles.  The symbols $a,a'$ will stand for generic colors.

The barbells
are the
broken words of the form
$$a^2\ \overset \leftarrow\ 2r + 1 \ \rightarrow \to  {\_ \ \ \dots \ \
\_}\ a^2$$
for  $0 \le r \le 2n-2$.

The racks are chosen from the broken words
of the form ~$HCT$ of length~$13n + 2J$, $1 \le J
\le V$, where  the {\it head\/}  is
$$H = (a\ \_)^I a^{2n - 2I}$$
for some~$I$, $1 \le I \le V$;
the {\it tail\/}  is
$$T = (\_\ a)^J$$
for some~$J$, $ 1 \le J \le V$;
and the {\it center\/}  is
$$C=a^{3n+i}\ \overset \leftarrow 2k \rightarrow \to
{\_ \dots\_}\ a\
\overset\leftarrow 2k \rightarrow \to {\_\dots \_}
\ a^{8n-4k-1-i}$$
for some $i$ and $k$, $0 \le i \le k-1$ and $1 \le k \le n$.

Given $I,J$ with $1 \le I,J \le V$, write
$$(A^{13n})_{IJ} = c_1(1!) + c_2(2!) + \dots + c_n(n!),$$
where $0 \le c_k \le k$ for $1 \le k \le n$.
If $c_k \ne 0$,
 choose the racks to be the $c_k = c_k(I,J)$ broken words  of the form
$$\left[(a\ \_)^I a^{2n - 2I}\right]\left[ a^{3n+i}\ \overset \leftarrow 2k
\rightarrow \to
{\_ \dots\_}\ a\
\overset\leftarrow 2k \rightarrow \to {\_\dots \_}
\ a^{8n-4k-1-i}\right]  \left[(\_\ a)^J\right]$$
for $0 \le i \le c_k-1$.

The barbells and racks have the following properties.

\noindent $\bullet$ The head $H = (a\ \_)^I a^{2n - 2I}$ of a rack
                    can be filled by the tail $T = (\_\ a')^J$
                    of a rack
                    in a tiling
                    if and only if $I=J$.

\noindent $\bullet$  Label the blanks in the center
 																				$$C = a^{3n+i}\overset\leftarrow 2k   \rightarrow \to
{\_ \dots\_}\ a\
                       \overset\leftarrow 2k \rightarrow \to {\_\dots \_}\
a^{8n-4k-1-i}$$
                      of a rack by $\{1,2,\dots,4k\}$.
	                     Barbells can appear in a tiling only in the gaps in
the centers
                      of  racks,
                      starting only in  odd places and  straddling the ~$a$.
                      Furthermore, the blanks in this center can be tiled
by barbells
                      in exactly $k!$~ways. To see this, define a permutation
                      ~$\pi$ of $\{1,2,\dots,k\}$ by $\pi(j) = \ell$ if and
only if
                       a barbell $a'a'\ \overset \leftarrow\ 2r + 1 \
\rightarrow
                       \to {\_ \ \ \dots \ \ \_}\ a'a'$ occupies places labelled
                       $2j-1,2j,2(k+\ell)-1,2(k+\ell)$.

\noindent $\bullet$ The heads of racks can appear in a tiling starting only at
                    places which differ by multiples of~$13n$.

\medskip

Let $T$ be the tiling system with   prototiles the barbells and racks chosen
above. Then $T =T_0 \cup T_1 \cup \dots \cup
T_{13n-1}$, where $T_i$ is the set of indexed bi-infinite sequences in~$T$
in which the heads
appear starting at places congruent to~$i$ modulo~$13n$.
Thus $T_0$ consists of all indexed bi-infinite concatenations of words of
length~$13n$, of
the form $\bar H \bar C$, starting at multiples of~$13n$, where $H$ and~$C$
are the head and center of a rack, and
$\bar H$ and $\bar C$ are the solid words resulting from
filling them in a tiling.
Recall that if $H = (a\ \_)^I a^{2n - 2I}$,
then it must be filled by a tail
$T = (\_\ a')^I$.
$C$ can be filled only by barbells.

Define an edge shift as follows.
Let $G'$ be the directed graph with
vertices $1,2,\dots,V$,
and an arc from~$I$ to~$J$ for each rack
with head $(a\ \_)^I a^{2n - 2I}$,
tail $(\_\ a)^J$,
and center tiled by barbells.
There are $(A^{13n})_{IJ}$ arcs from~$I$ to ~$J$,
and an arc with head $(a\ \_)^I a^{2n - 2I}$
can follow an arc with tail $(\_\ a')^J$
if and only if $I=J$.
Therefore since $m = 13n$, the adjacency matrix of~$G'$ is~$A^m$
and $\si^m : T_0 \to T_0$
is  topologically conjugate to
$\si : \Sigma_{A^m} \to \Sigma_{A^m}$,
which in turn is topologically conjugate to
$\si^m : \Sigma_A \to \Sigma_A$.
\qed
\enddemo

\proclaim{Corollary}
The set of topological entropies of \ts s   is the same as that of  shifts
of finite
type, i@.e@., the set of logarithms of roots of Perron numbers.
\endproclaim

\vfil\eject

\enddocument